\documentclass[11pt,aps,prb, preprint]{amsart}
\usepackage{amsmath}
\usepackage{amssymb}

\newtheorem{theorem}{Theorem}[section]
\newtheorem{lemma}[theorem]{Lemma}
\newtheorem{proposition}[theorem]{Proposition}
\newtheorem{corollary}[theorem]{Corollary}

\begin{document}
\title{A remark on the global well-posedness of a modified critical quasi-geostrophic equation}
\author{Kazuo Yamazaki}  
\date{}
\maketitle

\begin{abstract}
The $\beta$-generalized quasi-geostrophic equation is studied in the range of $\alpha \in (0, 1), \beta \in (1/2, 1), 1/2 < \alpha + \beta < 3/2$. When $\alpha \in (1/2, 1), \beta \in (1/2, 1)$ such that $1 \leq \alpha + \beta < 3/2$, using the method introduced in [12] and [9], we prove global regularity of the unique and analytic solution and when $\alpha \in (0, 1/2), \beta \in (1/2, 1)$ such that $1/2 < \alpha + \beta < 1$, that there exists a constant such that $\lVert \nabla\theta_{0}\rVert_{L^{\infty}}^{2-2\alpha-2\beta}\lVert\theta_{0}\rVert_{L^{\infty}}^{2\alpha + 2\beta - 1} \leq c_{\alpha,\beta}$ implies global regularity. 

\vspace{5mm}

\end{abstract}
\footnote{2000MSC : 35B65, 35Q35, 35Q86}
\footnote{Department of Mathematics, Oklahoma State University, 401 Mathematical Sciences, Stillwater, OK 74078, USA}

\section{Introduction}
The $\beta$-generalized quasi-geostrophic equation QG$_{\alpha, \beta}$ proposed in [10] in a two-dimensional torus $\mathbb{T}^{2}$ is defined as follows:

\begin{equation}
\begin{cases}
\partial_{t}\theta + u\cdot\nabla\theta + \nu\Lambda^{2\alpha}\theta = 0,\\
u = \nabla^{\bot}(-\triangle)^{-\beta}\theta = \Lambda^{1-2\beta}\mathcal{R}^{\bot}\theta, \hspace{5mm} \theta(x,0) = \theta_{0}(x)
\end{cases}
\end{equation}

where $\theta$ represents liquid temperature, $\nu$ $>$ 0 the dissipative coefficient which hereafter we assume to be one, $\mathcal{R}$ a Riesz transform and the operator $\Lambda$ has its Fourier symbol $\widehat{\Lambda f} = \lvert \xi\rvert \hat{f}$. The range of $\alpha$ and $\beta$ considered in [10] is $\alpha \in [0, 1/2)$ and $\beta \in [1/2, 1]$; here we consider $\alpha \in (0, 1)$ and $\beta \in (1/2, 1)$ such that $1/2 < \alpha + \beta < 3/2$. When $\alpha = 0$ and $\beta = 1$, the model describes the evolution of the vorticity of a two dimensional damped inviscid incompressible fluid. The case $\beta = 1/2$ is the dissipative quasi-geostrophic equation (QG) from the geostrophic study of rotating fluids and has been extensively studied recently, e.g. [2], [3], [4], [6], [13] and references found therein. When $\beta = 0, \alpha = 1$, we find the magneto-geostrophic equation studied in [7] to be a meaningful generalization of this endpoint case. 

In particular, when $\alpha = 0$ and $\beta = 1$, (1) becomes Euler equation in vorticity form, $\alpha = \beta = 1/2$ the critical QG and (1) at $\beta = 1-\alpha$ and $\alpha \in (0, 1/2)$ was originally introduced in [5] as the critical MQG interpolating in-between. The authors in [5] showed the global existence of smooth solutions with $L^{2}$ initial data using the method introduced in [2]. Subsequently in [15], the authors showed that for any initial data in $H^{m}, m > 2$, there exists a unique global solution to (1) with $\beta = 1-\alpha$ in the case $\alpha \in (0,1)$; similar result for other active scalars is attainable (cf. [17]).

The purpose of this paper is twofold. Firstly, applying the method introduced in [12] and [9] we show the global regularity of the unique and analytic solution to (1) at $\beta = 1-\alpha, \alpha \in (0, 1/2)$ with an improved initial regularity condition in $H^{1}$; in [14] the author obtained this result by a different method. Secondly, we generalize further by considering the whole range of $1/2 < \alpha + \beta < 3/2$, which may be considered as the supercritical if $1/2 < \alpha + \beta < 1$ and subcritical if $1 < \alpha + \beta < 3/2$ according to the $L^{\infty}$ maximum principle shown in [4]. In the supercritical case, the author in [11] showed the eventual regularization of solutions to (1). Moreover, using extended Besov space, the Corollary 1.6 of [3] showed global regularity of the unique solution in the case $\beta \in (1-\alpha, 1], \alpha \in (0, 1/2)$. Finally, in [16] the authors showed in particular the global regularity of the unique solution to (1) with $\beta \in (0, 1/2), \alpha \in (1/2, 1]$ such that $1 < \alpha + \beta < 3/2$. Now let $\lVert\cdot\rVert_{s}$ denote the norm of $H^{s}$ while $\lVert\cdot\rVert_{L^{p}}$ that of $L^{p}$. Our main results read: 

\begin{theorem}
Let $\beta \in (1/2, 1), \alpha \in (0, 1)$ such that $1/2 < \alpha + \beta < 3/2$. If $\theta_{0} \in H^{s}, s \geq 3-2\beta - 2\alpha$, then there exists $T = T(\theta_{0}) > 0$ such that a solution $\theta(x,t)$ of (1) satisfies

\begin{eqnarray*}
&&\theta(x,t) \in C([0, T], H^{s}) \cap L^{2}([0, T], H^{s+\alpha})\\
&&t^{n/2}\theta(x,t) \in C((0, T], H^{s+n\alpha}) \cap L^{\infty}([0, T], H^{s+n\alpha})
\end{eqnarray*}

for every $n > 0$. The solution $\theta(x,t)$ is unique in case $\beta \in (1/2, 1), \alpha \in (0, 1)$ such that $1/2 < \alpha + \beta \leq 1$ and $\beta \in (1/2, 1), \alpha \in (1/2, 1)$ such that $1 < \alpha + \beta < 3/2$ as well as analytic in spatial variable for any $t > 0$.
\end{theorem}

\begin{corollary}
Let $\beta \in (1/2, 1), \alpha = 1-\beta$. If $\theta_{0} \in H^{s}, s \geq 3-2\beta-2\alpha$, then there exists a unique analytic solution that remains smooth for all time.
\end{corollary}

\begin{theorem}
Let $\beta \in (1/2, 1), \alpha \in (1/2, 1)$ such that $1 < \alpha + \beta < 3/2$. If $\theta_{0} \in H^{s}, s \geq 3-2\beta-2\alpha$, then the unique and analytic solution remains smooth for all time.
\end{theorem}

The proof of the Corollary is Theorem 1.1 and the discovery of an appropriate modulus of continuity (MOC) in [15]. We stress that our range of $\alpha$  in Theorem 1.3 is different from that in [3] while the range of $\beta$ different from that in [16]. Next, we consider the supercritical regime:

\begin{theorem}
Let $\beta \in (1/2, 1), \alpha \in (0, 1/2)$ such that $1/2 < \alpha + \beta < 1$. If $\theta_{0}(x) \in H^{s}, s\geq 3-2\beta-2\alpha$, then there exists a constant $c_{\alpha, \beta}$ that depends on $\alpha$ and $\beta$ such that 

\begin{equation*}
\lVert \nabla\theta_{0}\rVert_{L^{\infty}}^{2-2\alpha-2\beta}\lVert\theta_{0}\rVert_{L^{\infty}}^{2\alpha+2\beta - 1} \leq c_{\alpha,\beta}
\end{equation*}

implies that the unique analytic solution remains smooth for all time. 
\end{theorem}

In the actual proof of extending local solution to global in time, we will rely on the periodicity of the solution; however, it is well-known that the work in [6] allows us to drop this condition. It is of much interest if the initial regularity may be extended to critical Besov space (cf. [1], [18]). Now in Section 2 we prove Theorem 1.1 and then in Section 3, Theorem 1.3 and 1.4. 

\section{Proof of Theorem 1.1}
\subsection{For $s > 3-2\beta-2\alpha$}

We employ Galerkin approximation with $(e^{2\pi ikx})_{k=-N}^{N}$ the first (2N + 1) eigenfunctions of Laplacian. We take $P^{N}$ the projection onto the (2N+1)-dimensional subspace spanned by these basis and work on

\begin{equation}
\partial_{t}\theta^{N} = -P^{N}(u^{N}\cdot\nabla\theta^{N}) - \Lambda^{2\alpha}\theta^{N}, \hspace{5mm} \theta^{N}(x, 0) = P^{N}\theta_{0}(x)
\end{equation}

On the Fourier side is, up to some constants,

\begin{equation}
\partial_{t}\widehat{\theta}^{N}(k, t) = - \sum_{l+m=k, \lvert l\rvert,\lvert m\rvert,\lvert k\rvert \leq N}<l, m^{\bot}>(\frac{1}{\lvert m\rvert^{2\beta}} - \frac{1}{\lvert l\rvert^{2\beta}})\widehat{\theta}^{N}(m)\widehat{\theta}^{N}(l) - \lvert k\rvert^{2\alpha}\widehat{\theta}^{N}(k)
\end{equation}

where $<l, m^{\bot}> = l_{1}m_{2} - l_{2}m_{1}$. We multiply (1) by $\Lambda^{2s}\theta^{N}$ and estimate

\begin{equation}
S := \sum_{l+m+k=0, \lvert l\rvert,\lvert m\rvert,\lvert k\rvert\leq N}<l, m^{\bot}>(\frac{1}{\lvert m\rvert^{2\beta}} - \frac{1}{\lvert l\rvert^{2\beta}})\lvert k\rvert^{2s}\hat{\theta^{N}}(k)\hat{\theta^{N}}(l)\hat{\theta^{N}}(m)
\end{equation}

We symmetrize over variables k, l and m to obtain

\begin{eqnarray*}
\lvert S\rvert&=&\frac{1}{3}\lvert\sum_{l+m+k=0, \lvert l\rvert,\lvert m\rvert, \lvert k\rvert \leq N}(<l, m^{\bot}>(\frac{1}{\lvert m\rvert^{2\beta}} - \frac{1}{\lvert l\rvert^{2\beta}})\lvert k\rvert^{2s}\\
&& +<m,k^{\bot}>(\frac{1}{\lvert k\rvert^{2\beta}} - \frac{1}{\lvert m\rvert^{2\beta}})\lvert l\rvert^{2s}\\
&& +<k, l^{\bot}>(\frac{1}{\lvert l\rvert^{2\beta}} - \frac{1}{\lvert k\rvert^{2\beta}})\lvert m\rvert^{2s})\widehat{\theta^{N}}(k)\widehat{\theta^{N}}(l)\widehat{\theta^{N}}(m)\rvert\\
&\leq& C\sum_{l+m+k=0, \lvert l\rvert\leq\lvert m\rvert\leq\lvert k\rvert}\lvert <l, m^{\bot}>(\frac{1}{\lvert m\rvert^{2\beta}} - \frac{1}{\lvert l\rvert^{2\beta}})\lvert k\rvert^{2s}\\
&+& <m, k^{\bot}>(\frac{1}{\lvert k\rvert^{2\beta}} - \frac{1}{\lvert m\rvert^{2\beta}})\lvert l\rvert^{2s}\\
&+& <k, l^{\bot}>(\frac{1}{\lvert l\rvert^{2\beta}} - \frac{1}{\lvert k\rvert^{2\beta}})\lvert m\rvert^{2s}\rvert\lvert \widehat{\theta^{N}}(k)\widehat{\theta^{N}}(l)\widehat{\theta^{N}}(m)\rvert
\end{eqnarray*}

Now we observe that $<m, k^{\bot}> = <k, l^{\bot}> = <l, m^{\bot}>$ and hence

\begin{eqnarray*}
&&\lvert <l, m^{\bot}>(\frac{1}{\lvert m\rvert^{2\beta}} - \frac{1}{\lvert l\rvert^{2\beta}})\lvert k\rvert^{2s}\\ 
&&+ <m, k^{\bot}>(\frac{1}{\lvert k\rvert^{2\beta}} - \frac{1}{\lvert m\rvert^{2\beta}})\lvert l\rvert^{2s} + <k, l^{\bot}>(\frac{1}{\lvert l\rvert^{2\beta}} - \frac{1}{\lvert k\rvert^{2\beta}})\lvert m\rvert^{2s}\rvert\\
&\leq& \lvert <l, m^{\bot}>\rvert\lvert \frac{\lvert k\rvert^{2s}}{\lvert m\rvert^{2\beta}} - \frac{\lvert m\rvert^{2s}}{\lvert k\rvert^{2\beta}} + \frac{1}{\lvert l\rvert^{2\beta}}(\lvert m\rvert^{2s} - \lvert k\rvert^{2s}) + \lvert l\rvert^{2s}(\frac{1}{\lvert k\rvert^{2\beta}} - \frac{1}{\lvert m\rvert^{2\beta}})\rvert
\end{eqnarray*}

We note that under the condition that $\lvert l\rvert \leq \lvert m\rvert \leq \lvert k\rvert$, $k = -l-m$ gives $\lvert k\rvert \leq \lvert l \rvert + \lvert m\rvert \leq 2\lvert m\rvert$ and estimate

\begin{equation}
<l, m^{\bot}> \leq C\lvert l\rvert\lvert m\rvert
\end{equation}

\begin{equation}
\lvert \frac{\lvert k\rvert^{2s}}{\lvert m\rvert^{2\beta}} - \frac{\lvert m\rvert^{2s}}{\lvert k\rvert^{2\beta}}\rvert \leq \frac{\lvert -l\rvert\lvert k\rvert^{2s+2\beta - 1}}{\lvert m\rvert^{2\beta}\lvert k\rvert^{2\beta}}\leq C\lvert l\rvert^{1-2\beta}\lvert m\rvert^{s-1}\lvert k\rvert^{s}
\end{equation}

Similarly,

\begin{equation}
\frac{1}{\lvert l\rvert^{2\beta}}(\lvert m\rvert^{2s} - \lvert k\rvert^{2s}) \leq \lvert l\rvert^{1-2\beta}\lvert k\rvert^{2s-1} \leq C\lvert l\rvert^{1-2\beta}\lvert m\rvert^{s-1}\lvert k\rvert^{s}
\end{equation}

Moreover,

\begin{equation*}
\lvert l\rvert^{2s}(\frac{1}{\lvert k\rvert^{2\beta}} - \frac{1}{\lvert m\rvert^{2\beta}}) \leq \frac{\lvert -l\rvert^{2s+1}\lvert k\rvert^{2\beta - 1}}{\lvert k\rvert^{2\beta}\lvert m\rvert^{2\beta}}\leq C\lvert m\rvert^{s-1}\lvert k\rvert^{s}\lvert l\rvert^{1-2\beta}
\end{equation*}

Combining this with (5), (6) and (7) together, we have

\begin{eqnarray}
\lvert S\rvert&\leq&C\sum_{l+m+k=0}\lvert k\rvert^{s}\lvert m\rvert^{s}\lvert l\rvert^{2-2\beta}\lvert\hat{\theta^{N}}(l)\hat{\theta^{N}}(m)\hat{\theta^{N}}(k)\rvert\\
&\leq& C\lVert\theta^{N}\rVert_{s+\delta}^{2}\sum_{l}\lvert l\rvert^{2 - 2\beta - 2\delta}\lvert\hat{\theta^{N}}(l)\rvert\nonumber\\
&\leq& C\lVert\theta^{N}\rVert_{s+\delta}^{2}(\sum_{l}(\lvert l\rvert^{q}\lvert\hat{\theta}^{N}(l)\rvert)^{2})^{1/2}(\sum_{l\neq 0}\lvert l\rvert^{2r})^{1/2}\leq C\lVert\theta^{N}\rVert_{s+\delta}^{2}\lVert\theta^{N}\rVert_{q}\nonumber
\end{eqnarray}

for $q > 3-2\beta-2\delta$ where we used Young's inequality for convolution and Parseval's formula. Thus, we have 

\begin{equation*}
\partial_{t}\lVert\theta^{N}\rVert_{s}^{2} \leq C\lVert\theta^{N}\rVert_{s+\alpha-\epsilon}^{2}\lVert\theta^{N}\rVert_{q} - 2\lVert\theta^{N}\rVert_{s+\alpha}^{2}
\end{equation*}

where we wrote $\delta = \alpha - \epsilon$ for $\epsilon > 0$ to be specified below. Now if $q \geq s + \alpha - \epsilon$, then immediately we have

\begin{equation*}
\partial_{t}\lVert\theta^{N}\rVert_{s}^{2} \leq C\lVert\theta^{N}\rVert_{q}^{M(q,s)} - \lVert\theta^{N}\rVert_{s+\alpha}^{2}
\end{equation*}

and if $q < s + \alpha - \epsilon$, then by interpolation

\begin{equation*}
\lVert \theta^{N}\rVert_{s+\alpha - \epsilon}^{2} \leq \lVert \theta^{N}\rVert_{s+\alpha}^{2(1-\gamma)}\lVert\theta^{N}\rVert_{q}^{2\gamma}
\end{equation*}

with $\gamma = \frac{\epsilon}{s+\alpha-q}$. Thus, with Young's inequality, 

\begin{equation*}
\partial_{t}\lVert\theta^{N}\rVert_{s}^{2} \leq C\lVert\theta^{N}\rVert_{q}^{2+\frac{1}{\gamma}} - (1+\gamma)\lVert\theta^{N}\rVert_{s+\alpha}^{2} \leq C\lVert\theta^{N}\rVert_{q}^{M(q,s)} - \lVert\theta^{N}\rVert_{s+\alpha}^{2}
\end{equation*}

Finally, if s = q, then taking $\gamma = \frac{\epsilon}{\alpha}$, the above interpolation gives

\begin{eqnarray*}
\partial_{t}\lVert\theta^{N}\rVert_{s}^{2} &\leq& C\lVert\theta^{N}\rVert_{s}(\lVert\theta^{N}\rVert_{s+\alpha}^{2(1-\frac{\epsilon}{\alpha})}\lVert\theta^{N}\rVert_{s}^{2(\frac{\epsilon}{\alpha})}) - 2\lVert\theta^{N}\rVert_{s+\alpha}^{2}\\
&\leq& C\lVert\theta^{N}\rVert_{s}^{2+\frac{\alpha}{\epsilon}} - \lVert\theta^{N}\rVert_{s+\alpha}^{2}
\end{eqnarray*}

Thus, we have shown, 

\begin{lemma}
For $q > 3-2\beta - 2\alpha, s \geq 0$ if $\theta_{0} \in H^{s}$, then

\begin{equation}
\partial_{t}\lVert\theta^{N}\rVert_{s}^{2} \leq C(q)\lVert \theta^{N}\rVert_{q}^{M(q,\alpha,s)} - \lVert\theta^{N}\rVert_{s+\alpha}^{2}
\end{equation}

and if s = q, then for $\epsilon \in (0, \min(\frac{q+2\beta + 2\alpha - 3}{2}, \alpha))$,

\begin{equation}
\partial_{t}\lVert\theta^{N}\rVert_{s}^{2} \leq C(\epsilon)\lVert\theta^{N}\rVert_{s}^{2+\frac{\alpha}{\epsilon}} - \lVert\theta^{N}\rVert_{s+\alpha}^{2}
\end{equation}
\end{lemma}

As a consequence of (10) and local existence of the solution to 

\begin{equation}
z' = Cz^{1+\frac{\alpha}{2\epsilon}}, \hspace{5mm} z(0) = z_{0}
\end{equation}

we have 

\begin{lemma}
For $s > 3-2\beta - 2\alpha$, if $\theta_{0} \in H^{s}$, there exists time $T = T(s, \alpha, \beta, \lVert\theta_{0}\rVert_{s}$) such that for every N uniformly we have the bound

\begin{equation*}
\lVert \theta^{N}\rVert_{s}(t) \leq C(s, \alpha, \beta, \lVert\theta_{0}\rVert_{s}), \hspace{5mm} 0 < t \leq T.
\end{equation*}
\end{lemma}

Next, we obtain uniform bounds for higher order of $H^{s}$ norms:

\begin{lemma}
Under the hypothesis of Lemma 2.2, there exists time $T = T(s, \alpha, \beta, \lVert\theta_{0}\rVert_{s})$ such that for all N uniformly 

\begin{equation}
t^{n/2}\lVert \theta^{N}\rVert_{s+n\alpha}\leq C(n, s, \alpha, \beta, \lVert\theta_{0}\rVert_{s}), 0 < t \leq T,
\end{equation}

for any $n \geq 0$. 
\end{lemma}

\proof{
We induct on n in integers and then interpolate. For n = 0, we see that it is done by Lemma 2.2. Now assume it is true for n; i.e.

\begin{equation}
\lVert\theta^{N}\rVert_{s+n\alpha}^{2}\leq Ct^{-n}
\end{equation}

Fix any $t \in [0, T]$ and consider an interval I$ = (t/2, t)$. By (9) we have 

\begin{equation}
\partial_{t}\lVert\theta^{N}\rVert_{s+n\alpha}^{2} \leq C(q)\lVert\theta^{N}\rVert_{s}^{M} - \lVert\theta^{N}\rVert_{s+(n+1)\alpha}
\end{equation}

and hence an integration in the interval I = (t/2, t) gives us 

\begin{eqnarray*}
\int_{I}\lVert\theta^{N}\rVert_{s+(n+1)\alpha}^{2}ds&\leq&C\int_{I}\lVert\theta^{N}\rVert_{s}^{M}ds + \lVert\theta^{N}(t/2)\rVert_{s+n\alpha}^{2}\\
&\leq&ct + \lVert\theta^{N}(t/2)\rVert_{s+n\alpha}^{2} \leq ct^{-n}
\end{eqnarray*}

where we used Lemma 2.2 and the induction hypothesis. Considering the average over I, we see that there exists some time $\eta$ in I such that 

\begin{equation}
\lVert\theta^{N}(\eta)\rVert_{s+(n+1)\alpha}^{2} \leq ct^{-n-1}
\end{equation}

Moreover, (14) gives 

\begin{equation*}
\partial_{t}\lVert\theta^{N}\rVert_{s+(n+1)\alpha}^{2} \leq C(q)\lVert\theta^{N}\rVert_{s}^{M}
\end{equation*}

and thus integration over $[\eta, t]$ gives us 

\begin{eqnarray*}
\lVert\theta^{N}(t)\rVert_{s+(n+1)\alpha}^{2}&\leq& C\int_{\eta}^{t}\lVert\theta^{N}\rVert_{s}^{M}dr + \lVert\theta^{N}(\eta)\rVert_{s+(n+1)\alpha}^{2}\\
&\leq& ct + \lVert\theta^{N}(\eta)\rVert_{s+(n+1)\alpha}^{2}\leq ct^{-n-1}
\end{eqnarray*}

Now for any $r \in \mathbb{R}^{+}, 0 < r \leq n$, Gagliardo-Nirenberg inequality completes the interpolation and the proof.}

\vspace{5mm}

Looking at (2) and (12), we see that for all $\epsilon > 0$ small and any $r > 0$, uniformly in N and $t \in [\epsilon, T]$, $\lVert\theta^{N}_{t}\rVert_{r}\leq C(r,\epsilon)$. With this and (12), the well-known compactness criteria implies that there exists a subsequence $\theta^{N_{j}}$ converging in $C([\epsilon, T], H^{r})$ to $\theta$. By the arbitrariness of $\epsilon$ and r, one can apply the standard subsequence of subsequence procedure to find a subsequence that converges to $\theta$ in C((0, T], $H^{r}$) for any r $> 0$. The limiting function $\theta$ still satisfies (12) and solves (1) on (0, T]. 

In order to show that $\theta$ converges to $\theta_{0}$ strongly in $H^{s}$ as $t \to 0$, we introduce $\phi(x)$ an arbitrary $C^{\infty}$ function and consider $g^{N}(t, \phi) \equiv (\theta^{N}, \phi) = \int\theta^{N}(x, t)\phi(x)dx$. Notice $g^{N}(\cdot, \phi) \in C([0, \tau])$ where $\tau \equiv T/2$ and taking an inner product of (2) with $\phi$ we obtain

\begin{eqnarray*}
\lvert\partial_{t}g^{N}(t,\phi)\rvert &\leq& C\lVert u^{N}\cdot\nabla\phi\rVert_{L^{2}}\lVert\theta^{N}\rVert_{L^{2}} + \lVert\theta^{N}\rVert_{L^{2}}\lVert\phi\rVert_{2\alpha}\\
&\leq& C\lVert \Lambda^{1-2\beta} \mathcal{R}^{\bot}\theta^{N}\rVert_{\dot{H}^{2\beta-1}}\lVert\nabla\phi\rVert_{\dot{H}^{2-2\beta}}\lVert\theta^{N}\rVert_{L^{2}} + \lVert\theta^{N}\rVert_{L^{2}}\lVert\phi\rVert_{2\alpha}\\
&\leq& C\lVert \theta^{N}\rVert_{L^{2}}^{2}\lVert\phi\rVert_{3-2\beta} + \lVert\theta^{N}\rVert_{L^{2}}\lVert\phi\rVert_{2\alpha}
\end{eqnarray*}

where we used the classical estimate that for every divergence-free f

\begin{equation}
s < 1, t < 1, s + t > -1 \Rightarrow \lVert f\cdot\nabla g\rVert_{\dot{H}^{s+t-1}}\leq c\lVert f\rVert_{\dot{H}^{s}}\lVert\nabla g\rVert_{\dot{H}^{t}}
\end{equation}

for some constant that depends on s and t and that if $\sigma > 0$, then $H^{\sigma} \subset \dot{H}^{\sigma}$ and $2\beta - 1, 2-2\beta > 0$. Finally, we also used the bound on Riesz transform in $L^{p}, p \in (1, \infty)$. Thus, for any $\delta > 0$,

\begin{equation*}
\int_{0}^{\tau}\lvert g_{t}^{N}\rvert^{1+\delta}dt \leq C(\int_{0}^{\tau}\lVert\theta^{N}\rVert_{L^{2}}^{2(1+\delta)} \lVert\phi\rVert_{3-2\beta}^{1+\delta}dt + \int_{0}^{\tau}\lVert\theta^{N}\rVert_{L^{2}}^{1+\delta}\lVert\phi\rVert_{2\alpha}^{1+\delta}dt)
\end{equation*}

By (10) we have $\lVert\theta^{N}\rVert_{L^{2}} \leq C$ on $[0, \tau]$ and thus $\lVert g_{t}^{N}(\cdot,\phi)\rVert_{L^{1+\delta}} \leq C(\phi)$ for any $\delta > 0$. Thus, we see that the sequence $g^{N}(\cdot, \phi)$ is compact in $C([0, \tau])$ and hence we can pick a subsequence $g^{N_{j}}(\cdot, \phi)$ converging uniformly to $g(\cdot, \phi) \in C([0,\tau])$. By choosing an appropriate subsequence we can assume $g(t, \phi) = \int\theta(x,t)\phi(x)dx$ for $t \in (0, \tau]$.

Next, we can choose a subsequence $(N_{j})_{j}$ such that $g^{N_{j}}(t, \phi)$ has a limit for any smooth function $\phi$ from a countable dense set in $H^{-s}$. Due to the uniform control over $\lVert\theta^{N_{j}}\rVert_{s}$ on $[0, \tau]$, we see that $g^{N_{j}}(t, \phi)$ converges uniformly on $[0, \tau]$ for every $\phi \in H^{-s}$. Note, for any $t > 0$,

\begin{equation*}
\lvert(\theta-\theta_{0},\phi)\rvert\leq\lvert (\theta-\theta^{N_{j}}, \phi)\rvert + \lvert (\theta^{N_{j}} - \theta_{0}^{N_{j}}, \phi)\rvert + \lvert(\theta_{0}^{N_{j}} - \theta_{0},\phi)\rvert
\end{equation*}

where the first and third tend to zero for $N_{j}$ sufficiently large while the second as t approaches zero for any fixed $N_{j}$. Thus, by definition $\theta(\cdot, t) \to \theta_{0}(\cdot)$ as t $\to 0$ weakly in $H^{s}$. This implies $\lVert\theta_{0}(\cdot)\rVert_{s} \leq \liminf_{t\to 0}\lVert\theta(\cdot,t)\rVert_{s}$. On the other hand, by (10) we see that for every N, $\lVert\theta^{N}\rVert_{s}^{2}(t)$ is always below the graph of the solution to (11); thus, $\lVert\theta_{0}\rVert_{s} \geq \limsup_{t \to 0}\lVert\theta\rVert_{s}(t)$. This completes the proof of existence of the solution with $\theta_{0}\in H^{s}, s > 3-2\beta - 2\alpha$.

\subsection{For $s \geq 3-2\beta-2\alpha$} 

\vspace{1mm}

The proof of extending the previous result to $s \geq 3-2\beta - 2\alpha$ is very similar to that in section 2.1; we provide a sketch of the proof for completeness. Denote a Hilbert space of periodic functions by

\begin{equation}
H^{s,\phi} = \{f \in L^{2} : \lVert f\rVert_{H^{s,\phi}}^{2} = \sum_{n}\lvert n\rvert^{2s}\phi(\lvert n\rvert)^{2} \lvert\widehat{f}(n)\rvert^{2} < \infty\}
\end{equation}

for $\phi: [0, \infty) \to [1, \infty)$ some unbounded increasing function and repeat the Galerkin approximation to estimate S of (4) with $\phi(\lvert k\rvert)^{2}$; i.e.

\begin{eqnarray*}
\lvert S\rvert&\leq& C\lvert\sum_{l+m+k=0, \lvert l\rvert\leq\lvert m\rvert\leq \lvert k\rvert}<l, m^{\bot}>(\frac{1}{\lvert m\rvert^{2\beta}} - \frac{1}{\lvert l\rvert^{2\beta}})\lvert k\rvert^{2s}\phi(\lvert k\rvert)^{2}\\
&& +<m,k^{\bot}>(\frac{1}{\lvert k\rvert^{2\beta}} - \frac{1}{\lvert m\rvert^{2\beta}})\lvert l\rvert^{2s}\phi(\lvert l\rvert)^{2}\\
&& +<k, l^{\bot}>(\frac{1}{\lvert l\rvert^{2\beta}} - \frac{1}{\lvert k\rvert^{2\beta}})\lvert m\rvert^{2s}\phi(\lvert m\rvert)^{2}\rvert\lvert\widehat{\theta^{N}}(l)\widehat{\theta^{N}}(m)\widehat{\theta^{N}}(k)\rvert
\end{eqnarray*}

A similar procedure as in section 2.1 leads to 

\begin{equation*}
\lvert S\rvert \leq C \sum_{l+m+k=0, \lvert l\rvert \leq \lvert m\rvert \leq \lvert k \rvert}\lvert m\rvert^{s+\alpha}\lvert k\rvert^{s+\alpha}\lvert l\rvert^{2-2\beta-2\alpha}\phi(\lvert m\rvert)\phi(\lvert k\rvert)\lvert\widehat{\theta^{N}}(l)\widehat{\theta^{N}}(m)\widehat{\theta^{N}}(l)\rvert
\end{equation*}

from which we can obtain

\begin{equation}
\partial_{t}\lVert\theta^{N}\rVert_{H^{s,\phi}}^{2}\leq(C\epsilon\lVert\theta^{N}\rVert_{H^{q,\phi}} - 1)\lVert\theta^{N}\rVert_{H^{s+\alpha,\phi}}^{2} + C(M(\epsilon))
\end{equation}

for $q \geq 3-2\beta-2\alpha$. Considering this differential inequality in comparison to those of Lemma 2.1, the same procedure we ran in the case of $s > 3-2\beta-2\alpha$ leads to the identical result for $s \geq 3-2\beta-2\alpha$ and hence existence of smooth solution with initial data in $H^{3-2\beta - 2\alpha}$ can be proven; we refer interested readers to [12] and [9] for details here.

\subsection{Uniqueness}
\subsubsection{Case $1/2 < \alpha + \beta \leq 1, \beta \in (1/2, 1), \alpha \in (0,1)$}
Suppose $\theta^{1}$ and $\theta^{2}$ both solve (1) with $u^{1}, u^{2}$, and $\theta_{0}^{1}, \theta_{0}^{2} \in H^{3-2\beta - 2\alpha}$ respectively. We let $\theta = \theta^{1} - \theta^{2}, u = u^{1} - u^{2}$ and observe that $\partial_{t}\theta = -u^{1}\cdot\nabla\theta - u\cdot\nabla\theta^{2} - \Lambda^{2\alpha}\theta$ and hence taking $L^{2}$ inner product, we obtain

\begin{eqnarray*}
\frac{1}{2}\partial_{t}\lVert\theta\rVert_{L^{2}}^{2}&\leq& \lVert\Lambda^{1-2\beta}\mathcal{R}^{\bot}\theta\rVert_{L^{\frac{2}{2-\alpha-2\beta}}}\lVert\nabla\theta^{2}\rVert_{L^{\frac{2}{2\beta + \alpha - 1}}}\lVert\theta\rVert_{L^{2}} - \lVert\Lambda^{\alpha}\theta\rVert_{L^{2}}^{2}\\
&\leq& C\lVert\theta\rVert_{\alpha}\lVert\theta^{2}\rVert_{3-2\beta - \alpha}\lVert\theta\rVert_{L^{2}} - \lVert\Lambda^{\alpha}\theta\rVert_{L^{2}}^{2}
\end{eqnarray*}

where we used Riesz potential inequality. This leads to

\begin{equation*}
\partial_{t}\lVert\theta\rVert_{L^{2}}^{2} \leq C\lVert\theta^{2}\rVert_{3-2\beta - \alpha}^{2} \lVert\theta\rVert_{L^{2}}^{2} - \lVert\Lambda^{\alpha}\theta\rVert_{L^{2}}^{2}
\end{equation*}

Since $\theta^{2} \in L^{2}([0, T], H^{s+\alpha})$, by Gronwall's inequality, $\lVert\theta(\cdot, t)\rVert_{L^{2}}^{2} \equiv 0$ for all $t \in [0, T]$. 

\subsubsection{Case $1 < \alpha + \beta < 3/2, \alpha \in (1/2, 1), \beta \in (1/2, 1)$} Let $\psi = -\Lambda^{-2\beta}\theta$, take the scalar product with the difference equation in the previous case and estimate

\begin{equation*}
\int u^{1}\cdot\nabla\psi\theta dx \leq \lVert u^{1}\cdot\nabla\psi\rVert_{\dot{H}^{\beta-\alpha}}\lVert\psi\rVert_{\dot{H}^{\alpha+\beta}} \leq C\lVert u^{1}\rVert_{\dot{H}^{-\alpha-\epsilon+2}}\lVert\nabla\psi\rVert_{\dot{H}^{\beta-1+\epsilon}}\lVert\psi\rVert_{\alpha+\beta}
\end{equation*}

where we used the Holder's inequality and (16) with $\epsilon \in (1/2, \alpha)$ such that $3/2 > \beta + \epsilon$. We continue the estimate above by 

\begin{eqnarray*}
&&C\lVert \theta^{1}\rVert_{\dot{H}^{3-2\beta-\alpha-\epsilon}}\lVert\psi\rVert_{\dot{H}^{\beta+\epsilon}}\lVert\psi\rVert_{\alpha+\beta}\\
&\leq& C\lVert\theta^{1}\rVert_{3-2\beta-\alpha-\epsilon}\lVert\psi\rVert_{\dot{H}^{\beta}}^{\frac{\alpha-\epsilon}{\alpha}}\lVert\psi\rVert_{\alpha+\beta}^{\frac{\epsilon}{\alpha}}\lVert\psi\rVert_{\alpha+\beta}\\
&\leq&C\lVert\theta^{1}\rVert_{3-2\beta-\alpha-\epsilon}^{\frac{2\alpha}{\alpha-\epsilon}}\lVert\psi\rVert_{\dot{H}^{\beta}}^{(\frac{\alpha-\epsilon}{\alpha})(\frac{2\alpha}{\alpha-\epsilon})}+ \lVert\psi\rVert_{\alpha+\beta}^{(\frac{\alpha+\epsilon}{\alpha})(\frac{2\alpha}{\alpha+\epsilon})}\\
&\leq& C\lVert\theta^{1}\rVert_{3-2\beta-2\alpha}^{\gamma}\lVert\theta^{1}\rVert_{3-2\beta-\alpha}^{2}\lVert\psi\rVert_{\dot{H}^{\beta}}^{2} + \lVert\psi\rVert_{\alpha+\beta}^{2}
\end{eqnarray*}

for $\gamma$ such that $(\frac{2\alpha}{\alpha-\epsilon})(3-2\beta-\alpha-\epsilon) - 2(3-2\beta-\alpha) = \gamma(3-2\beta-2\alpha)$. This implies 

\begin{equation*}
\partial_{t}\lVert\Lambda^{\beta}\psi\rVert_{L^{2}}^{2} \leq C\lVert\theta^{1}\rVert_{3-2\beta-2\alpha}^{\gamma}\lVert\theta^{1}\rVert_{3-2\beta-\alpha}^{2}\lVert\psi\rVert_{\dot{H}^{\beta}}^{2}
\end{equation*}

Since $\theta^{1} \in C([0, T], H^{s})\cap L^{2}([0, T], H^{s+\alpha})$, Gronwall's inequality implies the desired result.

\subsection{Analyticity}
The proof of showing that the global solution to (1) with the initial data $\theta_{0} \in H^{s}$ for $s \geq 3-2\beta-2\alpha, \beta \in (1/2, 1), \alpha \in (0, 1), 1/2 < \alpha + \beta < 3/2$ is analytic for all $t > 0$ is also similar to that in section 2.1; we sketch it for completeness. Considering the Galerkin approximation (3) again, we let $\xi_{k}^{N}(t) = \widehat{\theta^{N}}(k,t)e^{\frac{1}{2}\lvert k\rvert^{2\alpha}t}$ and $\gamma_{l,m,k}=\frac{1}{2}(\lvert l\rvert^{2\alpha} + \lvert m\rvert^{2\alpha} - \lvert k\rvert^{2\alpha})$. We multiply (1) by $e^{\frac{1}{2}\lvert k\rvert^{2\alpha}t}$ to obtain

\begin{equation*}
\partial_{t}\xi_{k}^{N}(t) = C\sum_{l+m=k, \lvert l\rvert,\lvert m\rvert,\lvert k\rvert \leq N}e^{-\gamma_{l,m,k}t}<l,m^{\bot}>(\frac{1}{\lvert m\rvert^{2\beta}} - \frac{1}{\lvert l\rvert^{2\beta}})\xi_{l}^{N}\xi_{m}^{N} - \frac{1}{2}\lvert k\rvert^{2\alpha}\xi_{k}^{N}
\end{equation*}

Now consider Y(t) = $\sum\lvert k\rvert^{6}\lvert\xi_{k}^{N}(t)\rvert^{2}$. We have

\begin{eqnarray*}
\frac{dY(t)}{dt}&=&C Re(\sum_{l+m+k=0, \lvert l\rvert,\lvert m\rvert,\lvert k\rvert\leq N}
<l,m^{\bot}>(\frac{1}{\lvert m\rvert^{2\beta}} - \frac{1}{\lvert l\rvert^{2\beta}})\lvert k\rvert^{6}\xi_{l}^{N}\xi_{m}^{N}\xi_{k}^{N})\\
&+&C Re(\sum_{l+m+k=0, \lvert l\rvert,\lvert m\rvert,\lvert k\rvert \leq N}(e^{-\gamma_{l,m,k}t}-1)<l,m^{\bot}>(\frac{1}{\lvert m\rvert^{2\beta}} - \frac{1}{\lvert l\rvert^{2\beta}}))\lvert k\rvert^{6}\xi_{l}^{N}\xi_{m}^{N}\xi_{k}^{N})\\
&-& \sum_{k}\lvert k\rvert^{6+2\alpha}\lvert\xi_{k}^{N}\rvert^{2} = I_{1} + I_{2} + I_{3}
\end{eqnarray*}

On $I_{1}$ symmetrizing over l, m and k as done in section 2.1 gives

\begin{eqnarray*}
I_{1}&\leq& C Re(\sum_{l+m+k=0, \lvert l\rvert\leq \lvert m\rvert\leq \lvert k\vert}\lvert l\rvert^{2-2\beta}\lvert m\rvert^{3}\lvert k\rvert^{3}\xi_{l}^{N}\xi_{m}^{N}\xi_{k}^{N}) \leq CY\sum_{l}\lvert l\rvert^{2-2\beta}\lvert \xi_{l}^{N}\rvert \leq CY^{3/2}
\end{eqnarray*} 

By Taylor expansion, $\lvert e^{-\gamma_{l,m,k}t}-1\rvert \leq \lvert \gamma_{l,m,k}\rvert t\leq \min\{\lvert l\rvert,\lvert m\rvert\}t$. Thus, similarly as in section 2.1, we estimate

\begin{eqnarray*}
&&\lvert I_{2}\rvert\leq C\lvert\sum_{l+m+k=0, \lvert l\rvert,\lvert m\rvert,\lvert k\rvert \leq N}\min\{\lvert l\rvert,\lvert m\rvert\}t<l, m^{\bot}>(\frac{1}{\lvert m\rvert^{2\beta}} - \frac{1}{\lvert l\rvert^{2\beta}})\lvert k\rvert^{6}\xi_{l}^{N}\xi_{m}^{N}\xi_{k}^{N})\\
&\leq& Ct\sum_{l+m+k=0}\lvert l\rvert^{3-2\alpha-2\beta}\lvert m\rvert^{3+\alpha}\lvert k\rvert^{3+\alpha}\lvert \xi_{l}^{N}\rvert\lvert \xi_{m}^{N}\rvert\lvert \xi_{k}^{N}\rvert \leq CtY^{1/2}\sum_{k}\lvert k\rvert^{6+2\alpha}\lvert\xi_{k}^{N}\rvert^{2}
\end{eqnarray*} 

Therefore, combining $I_{1}, I_{2}$ and $I_{3}$ gives

\begin{equation*}
\frac{dY(t)}{dt}\leq C_{1}Y^{3/2}+(C_{2}Y^{1/2}t-1)\sum_{k}\lvert k\rvert^{6+2\alpha}\lvert\xi_{k}^{N}
\rvert^{2}
\end{equation*}

Note $Y(0) = \lVert\theta_{0}\rVert_{3}^{2}$. This implies that for time interval small enough, we have an upper bound on Y uniformly in N. By the blow-up criterion below, we know that the $H^{s}$ norm for any s $>0$ of any solution to (1) is bounded uniformly. Thus, for all $t_{0} > 0$ uniformly in N and $t > t_{0}$, we can repeat the process above and have the bound on $\sum_{k}\lvert\widehat{\theta^{N}}(k,t)\rvert^{2}e^{\delta\lvert k\rvert}$ for some small $\delta = \delta(t_{0}, \theta_{0}) > 0$. By construction of $\theta$, we know it satisfies the same bound. 

\section{Proof of Theorems 1.3 and 1.4}

\subsection{Blow-up Criterion}

We state a blow-up criterion which, using a standard commmutator estimate (cf. [8]) can be readily proven:

\begin{lemma}
Suppose the solution to (1) $\theta(x,t)$ satisfies $\lVert\nabla\theta(\cdot, t)\rVert_{L^{\infty}} \leq C$ for all time $t \in [0, T]$. Then for every $s > 0$, there exists a constant C(s) such that $\lVert \theta(\cdot, t)\rVert_{s} \leq C(s)$ for all time $t \in [0, T]$.
\end{lemma}

\subsection{Subcritical Case}
We first focus on the range of $\alpha \in (1/2, 1), \beta \in (1/2, 1), 1 < \beta + \alpha < 3/2$. We denote by $\xi = \lvert x-y\rvert$ interchangeably upon convenience. We define a MOC to be a continuous, increasing concave function $\omega: [0, \infty) \mapsto [0,\infty)$ with $\omega(0) = 0$. We say $\theta$ has a MOC $\omega$ if $\lvert\theta(x)-\theta(y)\rvert\leq \omega(\lvert x-y\rvert) \forall x, y \in \mathbb{T}^{2}$. The blowup criterion above and the following result due to [13] makes it clear that in order to show global regularity of $\theta$, it suffices to show that $\theta$ has a MOC $\omega$ for all $t > 0$. 

\begin{proposition}
If $\omega$ is a MOC for $\theta(x, t): \mathbb{T}^{2}$ $\to \mathbb{R}$ for all $t > 0$, then $\lvert\nabla\theta\rvert(x) \leq \omega'(0)$ for all x $\in \mathbb{T}^{2}$
\end{proposition}

For this reason, we shall construct a MOC $\omega$ such that $\omega'(0) < \infty$. Next,

\begin{proposition}
Assume $\theta$ has a strict MOC satisfying $\omega'' (0 +) = - \infty$ for all $t < T$; i.e. for all $x, y \in \mathbb{T}^{2}, \lvert\theta(x, t) - \theta(y, t)\rvert < \omega(\lvert x-y\rvert)$, but not for $t > T$. Then, there exists x, y $\in \mathbb{T}^{2}$, x $\neq$ y such that $\theta$(x, T) - $\theta$(y, T) = $\omega$($\lvert x - y\rvert$).
\end{proposition}

Thus, the only scenario in which a MOC $\omega$ is lost is if there exists $ T > 0$ such that $\theta$ has the MOC $\omega$ for all $ t \in [0, T]$ and two distinct points x and y such that $\theta$(x, T) - $\theta$(y, T) = $\omega$($\lvert$x - y$\rvert$). We rule out this possibility by showing that in such case, $\frac{\partial}{\partial t}[\theta(x, t) - \theta(y, t)]\lvert_{t = T} <$ 0. Let us write 

\begin{equation*}
\frac{\partial}{\partial t}[\theta(x) - \theta(y)]\rvert_{t=T}
= - [( u \cdot \nabla \theta)(x) - (u \cdot \nabla \theta)(y)] - [(\Lambda^{2\alpha} \theta)(x) - (\Lambda^{2\alpha} \theta)(y)]\rvert_{t=T}
\end{equation*}

Our agenda now is to first estimate the Convection and Dissipation terms, to be specific find upper bounds that depend on $\omega$. Then we will construct the MOC $\omega$ explicitly that assures us that the sum of the two terms is negative to reach the desired result. We have the following estimate on the convection term due to originally [13] and later generalized in [15]:

\begin{proposition}
If $\theta$ has a MOC $\omega$, then $u = \Lambda^{1-2\beta}\mathcal{R}^{\bot}\theta$ for any $\beta \in (0, 1)$ has a MOC

\begin{equation*}
\Omega(\xi) = C_{1}(\int_{0}^{\xi}\frac{\omega(\eta)}{\eta^{2-2\beta}}d\eta + \xi\int_{\xi}^{\infty}\frac{\omega(\eta)}{\eta^{3-2\beta}}d\eta)
\end{equation*}

for some constant $C_{1}$ that depends on $\beta$.
\end{proposition}
With that in mind, using Proposition 3.3, the following is clear:

\begin{equation*}
u \cdot \nabla \theta(x) - u \cdot \nabla \theta(y) \leq \lim_{h \searrow 0}\frac{\omega(\xi + h \Omega(\xi)) - \omega(\xi)}{h} = \Omega(\xi)\omega'(\xi)
\end{equation*}

We also borrow the result below, originally from [13], generalized in [19]:

\begin{equation*}
C_{2}[\int_{0}^{\xi/2}\frac{\omega(\xi + 2\eta) + \omega(\xi - 2\eta) - 2\omega(\xi)}{\eta^{1+2\alpha}} d\eta + \int_{\xi/2}^{\infty}\frac{\omega(2\eta + \xi) - \omega(2\eta - \xi) - 2\omega(\xi)}{\eta^{1+2\alpha}}d\eta]
\end{equation*}

As discussed in e.g. [19], it suffices to find $\lambda > 0$ such that $\omega_{\lambda}(\xi)$ is a MOC of $\theta_{0}(x)$; in the critical case, $\omega$ must be unbounded but not in the subcritical regime. We define

\begin{equation*}
\begin{cases}
\omega(\xi) = \xi - \xi^{r} \hspace{20mm} \xi \leq \delta\\
\omega'(\xi) = \frac{\gamma}{\xi^{2\alpha+2\beta - 1}} \hspace{15mm} \xi > \delta
\end{cases}
\end{equation*}

for $r \in (1, 2)$. We see that $\omega$ is continuous and $\omega(0) = 0$. It can be readily checked that the first derivative is positive for $\delta$ sufficiently small and hence increasing. Clearly $\omega'(0) < \infty$; the second derivative if $\xi \leq \delta$ is negative. We also have $\omega''(\xi) = \gamma(1-2\alpha-2\beta)\xi^{-2\alpha-2\beta} < 0$ as $1-2\alpha-2\beta < 0$. Moreover, notice $\lim_{\xi\to 0_{+}} \omega''(\xi) = -\infty$ as $r < 2$. Finally, $\omega'(\delta_{+}) = \gamma\delta^{-(2\alpha+2\beta-1)} < 1-r\delta^{r-1} = \omega'(\delta_{-})$ if we take $\gamma$ small enough as $r > 1$. We consider two different cases now:

\textbf{Case: $\xi \leq \delta$}: Because we have $\frac{\omega(\xi)}{\xi} = 1-\xi^{r-1} \leq \omega'(0) = 1$,

\begin{equation*}
\int_{0}^{\xi}\frac{\omega(\eta)}{\eta^{2-2\beta}}d\eta \leq \int_{0}^{\xi}\eta^{2\beta-1} d\eta= \frac{\xi^{2\beta}}{2\beta}
\end{equation*}

Moreover,

\begin{equation*}
\int_{\xi}^{\delta}\frac{\omega(\eta)}{\eta^{3-2\beta}}d\eta = \int_{\xi}^{\delta}\eta^{2\beta-2} - \eta^{r-(3-2\beta)}d\eta\leq \frac{\delta^{2\beta-1}}{2\beta - 1} - \frac{\xi^{2\beta - 1}}{2\beta - 1} \leq \frac{\delta^{2\beta - 1}}{2\beta - 1}
\end{equation*}

Finally, 

\begin{equation*}
\int_{\delta}^{\infty}\frac{\omega(\eta)}{\eta^{3-2\beta}}d\eta = \frac{\omega(\delta)\delta^{2\beta - 2}}{2-2\beta} + \frac{\gamma}{(2-2\beta)(2\alpha)}\delta^{-2\alpha} \leq \frac{\delta^{2\beta - 1}}{2-2\beta} + \frac{\delta^{1-2\alpha}}{(2-2\beta)(2\alpha)}
\end{equation*}

Thus, the estimate from the convection term is 

\begin{equation*}
C_{1}[\frac{\xi^{2\beta}}{2\beta} + \xi[\frac{\delta^{2\beta - 1}}{2\beta - 1} + \frac{\delta^{2\beta - 1}}{2-2\beta} + \frac{\delta^{1-2\alpha}}{(2-2\beta)(2\alpha)}]]
\end{equation*}

To estimate dissipation term, note $\omega(\xi - 2\eta) \leq \omega(\xi) - 2\omega'(\xi)\eta + 4\omega''(\xi)\eta^{2}$ by Taylor expansion and hence using concavity

\begin{eqnarray*}
&&C_{2}\int_{0}^{\xi/2}\frac{\omega(\xi + 2\eta) + \omega(\xi - 2\eta) - 2\omega(\xi)}{\eta^{1+2\alpha}}d\eta \leq C_{2}\int_{0}^{\xi/2}\frac{4\omega''(\xi)\eta^{2}}{\eta^{1+2\alpha}}d\eta = -C_{2}\xi^{r-2\alpha}
\end{eqnarray*}

Thus,

\begin{equation*}
\xi[\frac{C_{1}\xi^{2\beta - 1}}{2\beta} + C_{1}[\frac{\delta^{2\beta - 1}}{2\beta-1} + \frac{\delta^{2\beta - 1}}{2-2\beta} + \frac{\delta^{1-2\alpha}}{(2-2\beta)(2\alpha)}] - C_{2}\xi^{r-2\alpha-1}]
\end{equation*}

Note $r - 2\alpha - 1 < 0$ as $r<2<1+2\alpha$ and $1-2\alpha > r-2\alpha-1$ since $2 > r$. Therefore, letting $\delta \to 0$ and hence forcing $\xi \to 0$, we achieve negativity.  

\textbf{Case $\xi \geq \delta$}: We now estimate

\begin{eqnarray*}
\int_{0}^{\xi}\frac{\omega(\eta)}{\eta^{2-2\beta}}d\eta \leq \omega(\xi) \int_{0}^{\xi}\frac{1}{\eta^{2-2\beta}}d\eta \leq \omega(\xi)\frac{\xi^{2\beta - 1}}{2\beta - 1}
\end{eqnarray*}

For the other integral, we integrate by parts and obtain

\begin{equation*}
\int_{\xi}^{\infty}\frac{\omega(\eta)}{\eta^{3-2\beta}}d\eta \leq \omega(\xi)\xi^{2\beta - 2}(\frac{1}{2-2\beta} + \frac{1}{(2-2\beta)(2\alpha)})
\end{equation*}

where we took $\gamma$ small enough so that 

\begin{equation*}
\gamma \leq \frac{1}{2}\delta^{2\beta+2\alpha-1} \leq \delta^{2\beta+2\alpha-1} - \delta^{2\beta+2\alpha-2+r} = \omega(\delta)\delta^{2\beta+2\alpha-2} \leq \omega(\xi)\xi^{2\beta+2\alpha-2}
\end{equation*}

Note above also implies 

\begin{equation}
\frac{2}{2^{2\alpha+2\beta}}\gamma \leq \omega(\xi)\xi^{2\beta+2\alpha}
\end{equation}

Thus, we now have the bound on the convection term:

\begin{equation*}
\Omega(\xi)\omega'(\xi)\leq C_{1}\frac{\omega(\xi)}{\xi^{2\alpha}}[\frac{\gamma}{2\beta-1} + \gamma(\frac{1}{2-2\beta} + \frac{1}{(2-2\beta)(2\alpha)})]
\end{equation*}

On the dissipation term, we have $\omega(2\eta + \xi) - \omega(2\eta - \xi) \leq \omega(2\xi)$ and

\begin{equation*}
\omega(2\xi) = \omega(\xi) + \gamma\int_{\xi}^{2\xi}\frac{1}{\eta^{2\alpha+2\beta-1}}d\eta \leq \omega(\xi) + \gamma(2\xi)^{-2\alpha-2\beta} \leq \frac{3}{2}\omega(\xi)
\end{equation*}

using (19). Thus,

\begin{eqnarray*}
&&C_{2}\int_{\xi/2}^{\infty}\frac{\omega(2\eta + \xi) - \omega(2\eta - \xi) - 2\omega(\xi)}{\eta^{1 + 2\alpha}}d\eta \leq -C_{2}\int_{\xi/2}^{\infty}\frac{\omega(\xi)}{\eta^{1+2\alpha}}d\eta \leq -C_{2}\omega(\xi)\xi^{-2\alpha}
\end{eqnarray*}

In sum, we have for $\gamma$ sufficiently small,

\begin{equation*}
C_{1}\frac{\omega(\xi)}{\xi^{2\alpha}}[\frac{\gamma}{2\beta - 1} + \gamma(\frac{1}{2-2\beta} + \frac{1}{(2-2\beta)(2\alpha)}) - C_{2}] < 0
\end{equation*}

\subsection{Supercritical Case}

We now consider $\alpha \in (0, 1/2)$ and $\beta \in (1/2, 1)$ such that $1/2 < \alpha + \beta < 1$. We define for $s \in (\alpha + \beta, 1)$ and $r \in (1, 1+2\alpha)$

\begin{equation*}
\begin{cases}
\omega(\xi) = \xi - \xi^{r} \hspace{10mm} \xi \in [0, \delta]\\
\omega'(\xi) = \frac{\gamma\delta^{s}}{\xi^{s}} \hspace{14mm} \xi > \delta
\end{cases}
\end{equation*}

Checking each requirement of MOC is similar to the previous case.

\textbf{Case $0 \leq \xi \leq \delta$}: Similarly to before, we have

\begin{eqnarray*}
&&\int_{0}^{\xi}\frac{\omega(\eta)}{\eta^{2-2\beta}}d\eta \leq \int_{0}^{\xi}\frac{\omega(\eta)}{\eta}d\eta \leq \xi\\
&&\int_{\xi}^{\delta}\frac{\omega(\eta)}{\eta^{3-2\beta}}d\eta \leq \int_{\xi}^{\delta}\eta^{2\beta - 2}d\eta\leq \frac{\delta^{2\beta - 1}}{2\beta - 1}
\end{eqnarray*}

On the other integral, by integration by parts,

\begin{equation*}
\int_{\delta}^{\infty}\frac{\omega(\eta)}{\eta^{3-2\beta}}d\eta \leq \frac{\delta^{-1+2\beta}}{2-2\beta} + \frac{\gamma\delta^{s}}{2-2\beta}\int_{\delta}^{\infty}\eta^{-2+2\beta-s}d\eta \leq \delta^{-1+2\beta}[\frac{s+2-2\beta}{(2-2\beta)(s+1-2\beta)}] 
\end{equation*}

Similar computation as before shows that

\begin{equation*}
C_{2}\int_{0}^{\xi/2}\frac{\omega(\xi + 2\eta) + \omega(\xi - 2\eta) - 2\omega(\xi)}{\eta^{1+2\alpha}}d\eta \leq C_{2}\int_{0}^{\xi/2}\frac{4\omega''(\xi)\eta^{2}}{\eta^{1+2\alpha}}d\eta = -C_{2}\xi^{r-2\alpha}
\end{equation*}

Combining these inequalities we let $\delta \to 0$ and attain

\begin{eqnarray*}
\xi[C_{1} + C_{1}[\frac{1}{2\beta - 1} + \frac{s+2-2\beta}{(2-2\beta)(s+1-2\beta)}]\delta^{2\beta - 1} - C_{2}\xi^{r-1-2\alpha}] < 0
\end{eqnarray*}

\textbf{Case $\xi > \delta$}: We compute

\begin{eqnarray*}
\int_{0}^{\xi}\frac{\omega(\eta)}{\eta^{2-2\beta}}d\eta \leq \int_{0}^{\delta}\eta^{-1+2\beta}d\eta + \omega(\xi)\int_{\delta}^{\xi}\eta^{-2+2\beta}d\eta \leq \frac{\delta^{2\beta}}{2\beta} + \frac{\omega(\xi)\xi^{2\beta - 1}}{2\beta - 1}
\end{eqnarray*}

Now $\omega(\xi) \geq \omega(\delta) = \delta - \delta^{r} \geq \delta^{2\beta}$ if $\delta$ is small. Therefore, we have 

\begin{equation*}
\int_{0}^{\xi}\frac{\omega(\eta)}{\eta^{2-2\beta}}d\eta \leq \frac{\omega(\xi)}{2\beta} + \frac{\omega(\xi)\xi^{2\beta - 1}}{2\beta - 1} = \omega(\xi)[\frac{1}{2\beta}+\frac{\xi^{2\beta - 1}}{2\beta - 1}]
\end{equation*}

On the other hand, by integration by parts,

\begin{equation*}
\int_{\xi}^{\infty}\frac{\omega(\eta)}{\eta^{3-2\beta}}d\eta = \xi^{2\beta - 2}[\frac{\omega(\xi)}{2-2\beta} + \frac{\gamma\xi^{1-s}}{(2-2\beta)}\frac{\delta^{s}}{(s+1-2\beta)}]
\end{equation*}

By definition $\omega(\xi) = \delta - \delta^{r} - \frac{\gamma\delta}{1-s} + \frac{\gamma}{1-s}\xi(\frac{\delta}{\xi})^{s}$. For $\delta$ small enough, we have $\frac{1}{2} \geq \delta^{r-1}$ and hence $1-\delta^{r-1} \geq \frac{1}{2}$; thus, $\delta - \delta^{r} \geq \frac{\delta}{2}$ so that

\begin{equation}
\omega(\xi) \geq \frac{\delta}{2} - \frac{\gamma\delta}{1-s} + \frac{\gamma\xi^{1-s}}{1-s}\delta^{s} = \delta[\frac{1}{2} - \frac{\gamma}{1-s}] + \frac{\gamma\xi^{1-s}}{1-s}\delta^{s} \geq \frac{\gamma\xi^{1-s}}{1-s}\delta^{s}
\end{equation}

if we take $\frac{1-s}{2} \geq \gamma$. Thus, we conclude 

\begin{equation*}
\int_{\xi}^{\infty}\frac{\omega(\eta)}{\eta^{3-2\beta}}d\eta \leq \xi^{2\beta - 2}\omega(\xi)[\frac{1}{2-2\beta} + \frac{1-s}{(2-2\beta)(s+1-2\beta)}]\\
\end{equation*}

With this we have the estimate on the convection term to be

\begin{eqnarray*}
\omega(\xi)\xi^{-2\alpha}\delta^{s}\gamma[C_{3}\xi^{2\alpha -s} + C_{4}\xi^{2\beta + 2\alpha - 1-s} + C_{5}\xi^{2\beta + 2\alpha - 1-s}]
\end{eqnarray*}

On the dissipation term, similarly as before, using (20) we obtain

\begin{equation*}
\omega(2\xi) \leq \omega(\xi) + 2(1-s)(\frac{\delta}{\xi})^{s}\gamma\xi \leq \omega(\xi) + 2^{1-s}(\frac{\delta}{\xi})^{s}\omega(\xi)(\frac{\xi}{\delta})^{s}(1-s) <  2\omega(\xi)
\end{equation*}

Thus, the contribution from dissipation can be bounded again similarly as before by $-C_{2}\omega(\xi)\xi^{-2\alpha}$. Hence,

\begin{eqnarray*}
&&\omega(\xi)\xi^{-2\alpha}\delta^{s}[\gamma[C_{3}\xi^{2\alpha -s} + C_{4}\xi^{2\beta + 2\alpha - 1-s} + C_{5}\xi^{2\beta + 2\alpha - 1-s}] - C_{2}\delta^{-s}]\\
&\leq&\omega(\xi)\xi^{-2\alpha}\delta^{s}[[C_{3}\gamma^{2\alpha -s+1} + C_{4}\gamma^{2\beta + 2\alpha - s} + C_{5}\gamma^{2\beta + 2\alpha - s}] - C_{2}\delta^{-s}]
\end{eqnarray*}

because $2\alpha - s < 0$ and $2\beta + 2\alpha - 1 - s < 0$. Now take $\gamma$ small enough and because $2\alpha - s + 1 > 0, 2\beta + 2\alpha - s > 0$, we have negativity. Q.E.D.

Finally, considering how small the initial data must be follows from definition of $\omega$; we sketch it for completeness. We have

\begin{equation*}
\omega_{\lambda}(\xi) =
\begin{cases}
\lambda^{2(\alpha+\beta-1)}[\lambda\xi - (\lambda\xi)^{r}] = \lambda^{2\alpha + 2\beta - 1}\xi - \lambda^{2(\alpha+\beta-1)+r}\xi^{r} \hspace{5mm} \xi\lambda \in [0,\delta]\\
\lambda^{2(\alpha+\beta-1)}[\frac{\gamma\delta^{s}(\lambda\xi)^{1-s}}{1-s} + \delta -\delta^{r} - \frac{\gamma\delta}{1-s}] \hspace{32mm} \xi\lambda > \delta
\end{cases}
\end{equation*}

Now for x, y such that $\lambda\lvert x-y\rvert = \lambda\xi \leq \delta$, we have $\lvert\theta_{0}(x) - \theta_{0}(y)\rvert\leq \omega_{\lambda}(\xi)$ if we set $\lambda^{2\alpha+2\beta - 1} = 2\lVert\nabla\theta_{0}\rVert_{L^{\infty}}$. For the case of $\lvert x-y\rvert > \frac{\delta}{\lambda}$, we have $\lvert\theta_{0}(x) - \theta_{0}(y)\rvert\leq 2\lVert\theta_{0}\rVert_{L^{\infty}}$ and therefore, $\omega_{\lambda}$ is a MOC of $\theta_{0}$ as long as 

\begin{equation*}
2\lVert\theta_{0}\rVert_{L^{\infty}} \leq \omega_{\lambda}(\frac{\delta}{\lambda}) = \lambda^{2(\alpha+\beta-1)}(\delta - \delta^{r}) = 2^{\frac{2(\alpha+\beta-1)}{2\alpha+2\beta-1}}\lVert\nabla\theta_{0}\rVert_{L^{\infty}}^{\frac{2(\alpha+\beta-1)}{2\alpha+2\beta-1}}(\delta-\delta^{r})
\end{equation*}

or equivalently $
\lVert\nabla\theta_{0}\rVert_{L^{\infty}}^{2(1-\alpha-\beta)}\lVert\theta_{0}\rVert_{L^{\infty}}^{2\alpha+2\beta-1}\leq 2^{-1}(\delta-\delta^{r})^{2\alpha+2\beta -1}$. 

\section{Acknowledgment}
The author expresses gratitude to Professor Jiahong Wu and Professor David Ullrich for their teaching and Professor Changxing Miao and Liutang Xue for their helpful comments.

\end{document}